\documentstyle[12pt]{article}
\begin{document}
\newtheorem{thm}{Theorem}
\newtheorem{prop}{Proposition}
\newtheorem{cor}{Corollary}
\newtheorem{ex}{Example}
\newtheorem{exs}{Examples}
\newtheorem{lem}{Lemma}
\newtheorem{rem}{Remark}
\newcommand{\bt}{\begin{thm}}
\newcommand{\et}{\end{thm}}
\newcommand{\bl}{\begin{lem}}
\newcommand{\el}{\end{lem}}
\newcommand{\bp}{\begin{prop}}
\newcommand{\ep}{\end{prop}}
\newcommand{\bc}{\begin{cor}}
\newcommand{\ec}{\end{cor}}
\newcommand{\p}{{{\bf Proof.\,\,}}}
\def\NN{I\!\! N}
\def\RR{I\!\! R}
\def\HH{I\!\! H}
\def\QQ{I\!\!\!\! Q}
\def\CC{I\!\!\!\! C}
\def\ZZ{Z\!\!\! Z}

\newcommand{\mod}[3]{\mbox{${#1} \equiv{#2}\bmod{#3}$}}
\newcommand{\G} {{\bf G}}
\newcommand{\biindice}[3]%
{
\renewcommand{\arraystretch}{0.5}
\begin{array}[t]{c}
#1\\
{\scriptstyle #2}\\
{\scriptstyle #3}
\end{array}
\renewcommand{\arraystretch}{1}
}
\newcommand{\chit}{\mbox{${\tilde\chi}$}}
\def\qed{\hfill \vrule height4pt width4pt depth2pt}
\baselineskip6mm \centerline{{\Large Strong commutativity preserving
maps  }} \vspace*{0.3cm}
 \centerline{{\Large on Lie ideals of semiprime rings }} \vspace*{0,6cm}
\centerline{{\large L. Oukhtite, S. Salhi \& L. Taoufiq}}
\centerline{\small{Universit\'e Moulay Isma\"\i l, Facult\'e des
Sciences et Techniques}} \centerline{\small{D\'epartement de
Math\'ematiques, Groupe d'Alg\`ebre et Applications}}
\centerline{\small{ B. P. 509 Boutalamine, Errachidia; Maroc}}
\centerline{\small{oukhtitel@hotmail.com, salhi@fastmail.fm,
lahcentaoufiq@yahoo.com}}
\begin{abstract}
\vspace*{-0.4cm} Let $R$ be a $2$-torsion free semiprime
ring and $U$ a nonzero square closed Lie ideal of $R$. In this paper
it is shown that if $f$ is either an endomorphism or an antihomomorphism
of $R$ such that $f(U)=U,$ then $f$ is strong commutativity preserving on $U$ if and
only if $f$ is centralizing on $U.$

\end{abstract}
semiprime rings,
Lie ideals.\\

 \section {Introduction}
 Throughout the present paper $R$ will denote a unitary associative
 ring. As usual, for
 $x, y $ in $R$, we write $[x,y]=xy-yx $, and we will use the
 identities $[xy,z]=x[y,z]+[x,z]y$, $[x,yz]=[x,y]z+y[x,z]$. For any $a\in R$, $d_a$ will
 denote the inner-derivation defined by $d_a(x)=[a,x]$ for all $x\in R$.\\
A ring  $R$  is said to be semiprime if $aRa=0$ implies that a=0. An
ideal $P$ of $R$ is prime if $aRb\subseteq P$ implies that $a\in P$
or $b\in P$. Recall that a ring $R$ is semiprime if and only if  its zero ideal
is the intersection of its prime ideals. Moreover, if the zero ideal
of $R$ is prime, then $R$ is said to be a prime ring.  An additive subgroup $U$ of a ring $R$ is a Lie ideal if
$[U,R]\subseteq U.$ Moreover, if $u^2\in U$ for all $u\in U$, then
$U$ is called a square closed Lie ideal. Since $(u+v)^2\in U$ and
[$u,v]\in U,$ we see that $2uv\in U$ for all $u, v\in U$. For a
subset $S$ of $R$, denote by $ann_R(S)$ the two-sided annihilator
 of $S$ -i.e.  $\{x\in R / Sx=xS=\{0\}\}.$  For every ideal $J$ of a semiprime ring $R$, it is known that $ann_R(J)$ is invariant under all derivations and  $J\cap
 ann_R(J)=0.$\\
 A  map $f:R \longrightarrow R$ is centralizing
 on $S$ if $[f(x),x]\in Z(R)$ for all $x\in S$; in particular if
 $[f(x),x]=0$ for all $x\in S$, then $f$ is called commuting on
 $S$. \\
A  map $f: R \longrightarrow R$ is called commutativity preserving
on $S$ if $[f(x),f(y)]=0$ whenever $[x,y]=0$, for all $x,y\in S.$ In
particular, if $[f(x),f(y)]=[x,y]$ for all $x,y\in S,$ then $f$ is
called strong commutativity preserving on $S.$\\Recently, M. S.
 Samman  \cite{MS} proved that an epimorphism
 of a semiprime ring is strong commutativity preserving if and only
 if it is centralizing on the entire ring. Moreover, he proved that if $R$ is a $2$-torsion free
 semiprime ring, then a centralizing antihomomorphism of $R$ onto itself must be strong commutativity preserving.
 The purpose of this paper is to extend the results of  \cite{MS} to
 square closed Lie ideals.
 \vspace*{-0.5cm}
 \section{Preliminaries and results}
 In order to prove our main theorems, we shall need the following results.
 \vspace*{-0.4cm}
  \bl \label{l1}

  Let $R$ be a $2$-torsion free semiprime ring and $U$ a nonzero Lie ideal of $R.$
If $[U,U]=0$, then  $U\subseteq Z(R).$

 \el
 \vspace*{-0.5cm}
\p
Let $u\in U$; since $[u,rt]\in U$ for all $r,t\in R$, then
  $[u,[u,rt]]=0$. Hence $u[u,rt]=[u,rt]u$. Therefore
  $$ur[u,t]+u[u,r]t=r[u,t]u+[u,r]tu.$$
  As $u[u,r]=[u,r]u$ and $[u,t]u=u[u,t]$, then
  $$ ur[u,t]+[u,r]ut=ru[u,t]+[u,r]tu.$$
  It follows that $2[u,r][u,t]=0$ for all $u\in U$ and $r, t \in R$.
 Since $R$ is 2-torsion-free,  thus
  \begin{equation} \label{equ1}
   [u,r][u,t]=0,\;\;\mbox{for all}\;\;u\in U\;\;\mbox{ and}\;\; r, t \in R.
\end{equation}
Replace $t$ by $sr$ in (1) to get
    $[u,r]R[u,r]=0$ for all $u\in U, r, t\in R$. The fact $R$ is
     semiprime implies that $U\subseteq Z(R)$.
     \qed\\
\\
    In all that follows $U$ will be a square closed Lie ideal of $R$
    and $M$ will denote the ideal of $R$ generated by $[U,U],$ that
    is $M=R[U,U]R.$
\vspace*{-0.4cm} \bl

Let $R$ be a $2$-torsion free semiprime ring and $d$ a derivation of
$R.$ If $a$ in $R$ satisfies $ad(U)=0$, then $ad(M)=0.$

\el
\p
Let $P$ be an arbitrary  prime ideal of $R,$ and note that
$\overline{R}=\displaystyle\frac{R}{P}$ is prime.  If $[U,U]\subseteq P$ or $char(\overline{R})=2,$ then
$2ad(R)M\subseteq P$ and $2Mad(R)\subset P.$ Assume now that $[U,U]\not\subset P$ and
$char(\overline{R})\neq 2.$ The fact that $R$ is 2-torsion-free and $ad(U)=\{0\}$ implies that $aUd(v)=\{0\}$ for all $v\in U$ and thus  $\bar{a}\overline{U}\overline{d(U)}=\bar{0}$. As
$[U,U]\not \subset P$, then $ \overline{U}\not\subset
Z(\bar{R}).$ Since  $[\overline{U}, \overline{U}]\neq \bar{0}$  from [4, Lemma 4]  either
$\overline{d(U)}=\bar{0}$ or $\bar{a}=\bar{0}$, that
is $d(U)\subseteq P$ or $a\in P .$ If $d(U)\subseteq P,$ then
$d[r,u]\in P$ for all $r\in R$ and $u\in U$. Replace $r$ by $rv$,
where $v\in U,$ to get $d(R)[U,U]\subseteq P$. Thus
$d(R)R[U,U]\subseteq P$ which yields $d(R)\subseteq P$ because
$[U,U]\not \subset P.$ In conclusion $ad(R)\subseteq P$.
Consequently, $ad(R)M\subseteq P$ and $Mad(R)\subseteq P.$ We now know that $2ad(R)M\subseteq P$ and $2Mad(R)\subseteq P$ for all prime ideals $P$ of $R,$ hence  $2ad(R)M=2Mad(R)=\{0\}.$ By 2-torsion-freeness we conclude that  $ad(R)M=Mad(R)=\{0\}.$ If we set  $J=ann_R(ann_R(M))$, then obviously $ad(R)J=0$. Since $R$ is semiprime, then $d(J)\subseteq J$ so that $ad(J)\subseteq J\bigcap
ann_R(J).$ Once again using the semiprimeness of $R$, we conclude
that $J\bigcap ann_R(J)=0$ so that $ad(J)=0.$ Since $M\subseteq J$,
this leads us to $ad(M)=0$. \qed
\vspace*{-0.5cm}
\bl \label{l3}

Let $R$ be a 2-torsion free semiprime ring. If $z\in U$ is such that
$z[U,U]=0$, then $[z,U]=0$.

 \el
  \vspace*{-0.5cm}
 \p
 If $[U,U]=0$, then $U\subseteq Z(R)$ by Lemma 1 and therefore $[z,U]=0.$
 Now suppose that $[U,U]\neq 0;$ from $z[U,U]=0$ we get
 $zd_u(v)=0$ for all $u,v\in U.$ Using Lemma 2,
 we find that $\;zd_u(x)=0\;$ for all $\;u\in U,\;x\in M=R[U,U]R.$ But
 $zd_u(x)=0$ assures that $zd_x(u)=0$ for all $u\in U, x\in M$ and once
 again using Lemma 2, we get $\;zd_x(M)=0,$ for all $\;x\in
 M.$ Hence $zd_x(y)=0$ for all $x,y\in M$ and thus   \vspace*{-0.3cm}$$z[x,y]=0\;\;\mbox{for
 all}\;\;x,y\in M.  \vspace*{-0.3cm}$$ Replace $y$ by $yz$ to get $zy[x,z]=0$, so that  $zM[x,z]=0
$. In view of $zxM[x,z]=0$, we then obtain $[x,z]M[x,z]=0.$ Since an ideal of a semiprime ring is semiprime,  $[x,z]=0$ for all $\;x\in
 M.$ As $R[U,U]\subseteq M,$ then $\;[z, r[u,v]]=0$ for all $r\in R,
 u,v\in U$. Using $[u,v]\in M,$ it then follows that
 $\;[z,r][u,v]=0.$ Replace $r$ by $rs$ in the least equality, we find
 that $[z,r]s[u,v]=0$ so that $\;[z,r]R[u,v]=0,$ for
 all $\;u,v\in U, r\in R.$ In particular $\;[z,v]R[z,v]=0$, proving
 $[z,v]=0$ for all $v\in U $ and  thus $[z,U]=0.$ \qed

\vspace*{-0.1cm}
\noindent  Now we are ready for our first theorem.
  \vspace*{-0.4cm}
\bt

Let $R$ be a 2-torsion free semiprime ring and $U$ a nonzero
 square closed Lie ideal of $R$. Suppose that $f$ is an endomorphism
 of $R$ such that $f(U)=U$. Then  $f$ is strong commutativity preserving
 on $U$ if and only if   $f$ is centralizing on $U$.

 \et
\p
From $[x,2xy]=[f(x),f(2xy)]$ for
all $x,y\in U,$ it follows that $(x-f(x))[x,y]=0$ for all $x,y\in
U$. Replacing $y$ by $2uy$ where $u, y\in U$, we get
\begin{equation} \label{equ1}
   (x-f(x))U[x,y]=0\;\;\mbox{for all}\;\;x,u\in U.
\end{equation}
As $2[U,U]R\subseteq U$ (because $2[u,v]r = 2[u,vr]- 2v[u,r]$), then
$(2)$ implies that
\begin{equation} \label{equ1}
  (x-f(x))[U,U]R[x,y]=0\;\;\mbox{for all}\;\;x,y\in U.
\end{equation} Let $P$ be an arbitrary prime ideal of $R$. It follows from (3) that for each $x\in U,$ either $(x-f(x))[U,U]\subseteq P$ or $[x, U]\subseteq P.$ The two sets of elements of $U$ for which these conditions hold are additive subgroups of $U$ whose union is $U$, hence  one must be equal to $U.$ Therefore $(x-f(x))[U,U]\subseteq P$ for all $x\in U$ and all prime ideals $P$-i.e., $(x-f(x))[U,U]=\{0\}$ for all $x\in U.$  Since $f(U)\subseteq U,$ then $u-f(u)\in U$  for all $u\in U$ and
Lemma 3 yields
$$[u-f(u),v]=0\;\;\mbox{for all}\;\;u,v\in U.$$Consequently, $\;\;[f(u),u]=0\;$ for all $\;u\in U\;$ so that $f$ is
commuting on $U.$ Accordingly, $f$ is centralizing on
$U.$\\Conversely, suppose that $[f(x),x]\in Z(R)$ for all $x\in U$.
By linearization $[x, f (y)]+[y, f(x)]\in  Z(R) $ for all $x,y$ in
$U.$ Using $[x, f (x^2)] + [x^2, f (x)]\in Z(R)$ together with
2-torsion-freeness, we find that $(x + f(x))[x, f (x)]\in Z(R),$ for
all $x\in U.$ Hence $[(x + f (x))[x, f (x)], x]=0$ and therefore
$[x, f (x)]^2=0.$ Since $[x, f (x)]$ in $Z(R)$, this yields
$[x, f (x)]R[x, f(x)] = 0$ and the semiprimeness of R forces $$[x, f
(x)] = 0\;\;\mbox{for all}\;\;x\in U.$$ Thus  $f$ is commuting on
$U$ and therefore $[f(x),y]=[x,f(y)]$ for all $x, y\in U.$ As $R$ is 2-torsion-free, then $
[f(x),xy]=[x,f(xy)] $ and thereby
$(f(x)-x)[f(x),y]=0$ for all $x, y\in U$. Replacing $y$ by $2uy$
where $u\in U$, we get $(f(x)-x)u[f(x),y]=0,$ so that
$(f(x)-x)U[x,f(y)]=0.$ Since $f(U)=U,$ then  $\;(f(x)-x)U[x,y]=0$
for all $x,y\in U.$ From $2[U,U]R\subseteq U$, it then follows that
$$(f(x)-x)[U,U]R[x,y]=0\;\;\mbox{for all}\;\;x,y\in U.$$
Reasoning as in the first part of the proof, we find that
$[f(z)-z,u]=0$ for all $z,u\in U,$ and therefore
$\;\;[f(z),u]=[z,u],$ for all $\;z, u\in U.$ Consequently, for $y, z
\in U$, this leads us to  $[f(z),f(y)]=[z,f(y)]=[z,y],$ proving that
$f$ is strong commutativity preserving on $U$.\qed
\\
\\
{\bf Remark.} From the proof of Theorem 1, one can easily see that
the condition $f(U)\subseteq U$ is sufficient to prove that $f$ is
strong commutativity preserving implies that $f$ is commuting on $U$
and therefore centralizing on $U.$\\
\\We easily derive the Proposition 2.1 of \cite{MS}, for $2$-torsion
free semiprime rings, as a corollary to Theorem 1.
 \vspace*{-0.4cm}
\bc

Let $f$ be an epimorphism of a $2$-torsion free semiprime ring $R.$
Then $f$ is strong commutativity preserving if and only if $f$ is
centralizing.

\ec
 \vspace*{-0.4cm}
 \noindent
In \cite{MY} it is proved that if $R$ is a $2$-torsion free prime
ring and $T$ an automorphism of $R$ which is centralizing on a Lie
ideal $U$ of $R$ and nontrivial on $U,$ then $U$ is contained in the
center of $R.$ Accordingly, in the special case when $U=R$, Theorem
2 gives a commutativity criterion as follows.
 \vspace*{-0.4cm}
\bc

Let $f$ be a nontrivial automorphism of a 2-torsion free prime ring
$R.$ If $f$ is strong commutativity preserving, then $R$ is
commutative.

\ec
  \vspace*{-0.3cm}
\noindent
To end this paper, the following theorem gives a condition
under which an antihomomorphism becomes strong commutativity
preserving. \vspace*{-0.4cm}

 \bt

 Let $R$ be a 2-torsion free semiprime ring and $U$ a
 square closed Lie ideal of $R$. If $f$ is an antihomomorphism
 of $R$ such that $f(U)= U$, then $f$ is centralizing on $U$ if and only if $f$
 is strong commutativity preserving  on $U.$

 \et
  \vspace*{-0.5cm}
\p Suppose $[U,U]\neq 0$ and then $M=R[U,U]R$ is a nonzero ideal of
$R.$ If $f$ is centralizing on $U$, then reasoning as in the proof
of Theorem 1 we find that $f$ is commuting on $U,$ so that
$[f(x),y]=[x,f(y)]$ for all $x, y\in U.$ Since $R$ is 2-torsion-free,  using
$[f(x), 2xy]=[x, f(2xy)]$ together with $f(U)=U$ we get
\begin{equation}
x[x,y]=[x,y]f(x)\;\;\mbox{for all}\;\;x,y\in U.
\end{equation}
Replace $y$ by $2uy$ in $(4),$ where $u\in U$, and once again using
2-torsion-freeness,  we get $\;[x,u][x,y+f(y)]=0.$ Write $2uv$ instead of
$u$ in this equality, with $v\in U,$ to find that
$\;[x,u]v[x,y+f(y)]=0.$  Hence
\begin{equation}
[x,u]U[x,y+f(y)]=0\;\;\mbox{for all}\;\;x,u, y\in U.
\end{equation}
Since $f(U)\subseteq U,$ replacing $u$ by $y+f(y)$ in $(5),  $ we
conclude that
\begin{equation}
[x,y+f(y)]U[x,y+f(y)]=0\;\;\mbox{for all}\;\;x,y\in U.
\end{equation}
If we set $T(U)=\{x\in R/ [x,R]\subseteq U\},$ then $[T(U),
R]\subseteq U\subseteq T(U)$ and from [\cite{Topic}, Lemma 1.4, p. 5] it follows that
$T(U)$ is a subring of $R.$ Moreover, $R[T(U), T(U)]R\subseteq
T(U).$ Indeed, let $x,y\in T(U)$ and $r\in R.$ From
$[x,yr]=[x,y]r+y[x,r]\in T(U)$ and $y[x,r]\in T(U)$ it follows
that $[x,y]r\in T(U).$ Since $[T(U), R]\subseteq T(U)$, then
$$[[x,y]r,s]=[x,y]rs-s[x,y]r \in T(U)\;\;\mbox{for all}\;\;r,s\in
R;$$ and therefore $\;s[x,y]r\in T(U)$ so that $R[T(U),
T(U)]R\subseteq T(U).$ In particular
$R[U,U]R\subseteq T(U)$, which proves that $\;[M,R]\subseteq U,$ where $M=R[U,U]R.$\\
In view of $(6)$, if we set $[x,y+f(y)]=a$ then $aUa=0.$ Let $u\in
U,\; m\in M$ and $r\in R$; from $[mau,r]\in [M,R]\subseteq U$ it
follows that $$
  0 = a[mau,r]a=a[ma,r]ua+ama[u,r]a=a[ma,r]ua=
   amarua,$$
so that $amaRua=0.$ Using $2am\in 2[U,U]R\subseteq U,$ we gLemma 1.4,et
$amaRama=0,$ hence  $\;aMa=0.$ Since $a\in M,$ we obviously get
$a=0,$ which implies that $[f(x),y]=[y,x],$ for all $\;x,y\in U.$
Accordingly,
$$[f(x),f(y)]=[f(y),x]=[x,y]\;\;\mbox{for all}\;\;x,y\in U,$$
proving that $f$ is strong commutativity
 preserving  on $U.$\\Conversely, if $f$ is strong commutativity preserving on
 $U,$ then \begin{equation}[f(x), f(y)]=[x,y],\;\;
 \mbox{for all}\;\;x,y\in U.\end{equation}
 Replace $y$ by $2xy$ in (7) we obtain
 \begin{equation}x[x, y]=[x, y]f(x).\end{equation}
Write $2uy$ instead of $y$ in (8), where $u\in U$, to find that
 $$xu[x,y] + x[x,u]y =u[x,y]f(x) + [x,u]yf(x).
$$
Since $x[x,u]y =[x,u]f(x)y$ and $[x,y]f(x)=x[x,y]$, by (8), then
$$xu[x,y] + [x,u]f(x)y =ux[x,y] + [x,u]yf(x)$$ and therefore
\begin{equation} [x,u][x+f(x), y]=0\;\;
 \mbox{for all}\;\;x,y,u \in U.\end{equation}
Replacing $y$ by $x$ in (9), we obtain
\begin{equation}[x,u][x,
f(x)]=0\;\;
 \mbox{for all}\;\;x,u \in U.
\end{equation} As $f(U)\subseteq U,$ write $2f(x)u$ instead of $u$ in
(10) to get $\;[x,f(x)]u[x,f(x)]=0\;$ and thus
$$
[x,f(x)]U[x, f(x)]=0.
$$
If we set $a=[x, f(x)],$ then $aUa=0$ and $a\in M=R[U,U]R.$
Reasoning as in the first part of our proof, we conclude that $a=0$
so that $[x, f(x)]=0$. Accordingly,  $f$ is commuting on $U$ and
therefore $f$ is centralizing on $U.$ \qed
\\
\\
{\bf Remark}. In the particular case when $U=R,$ the implication
that $f$ is strong commuativity preserving implying that $f$ is
centralizing is still valid without conditions on characteristic of
$R.$\\
\\
In [\cite{MS}, Proposition 2.4] M. S. Samman proved that if $R$ is a
$2$-torsion free semiprime ring, then a centralizing
antihomomorphism of $R$ onto itself must be
 strong commutativity preserving. Applying Theorem 2, we obtain a
 more general result as follows
 \vspace*{-0.2cm}
 \bc

Let $R$ be a 2-torsion free semiprime ring. If $f$ is an
antihomomorphism of $R$ onto itself, then $f$ is centralizing if and
only if $f$ is strong commutativity preserving.

 \ec

 \end{document}